\documentclass[11pt,a4paper]{article}

\usepackage{graphicx}
\usepackage{latexsym}
\usepackage{eepic}
\usepackage{amsfonts}

\makeatletter
\newdimen\mynormalparindent
\def\mymakefnmark{}
\def\mymakefntext{\indent\mymakefnmark}
\long\def\myfootnotetext#1{\insert\footins{%
  \normalfont\footnotesize
  \interlinepenalty\interfootnotelinepenalty
  \splittopskip\footnotesep \splitmaxdepth \dp\strutbox
  \floatingpenalty\@MM \hsize\columnwidth
  \@parboxrestore \parindent\mynormalparindent \sloppy
  \mymakefntext{\rule\z@\footnotesep\ignorespaces#1\unskip\strut\par}}}

\makeatother

\long\def\symbolfootnote[#1]#2{\begingroup%
\def\thefootnote{\fnsymbol{footnote}}\footnote[#1]{#2}\endgroup}

\newtheorem{basic}{Basic}[section]

\newtheorem{lem}[basic]{Lemma}

\newtheorem{thm}[basic]{Theorem}

\newcommand{\bdm}{\begin{displaymath}}
\newcommand{\edm}{\end{displaymath}}
\newcommand{\be}{\begin{equation}}
\newcommand{\ee}{\end{equation}}
\newcommand{\ep}{\vspace{-3mm}\hfill\mbox{$\Box$}\\}

\newcommand{\R}{\mathbb{R}}

\begin{document}

\begin{center}
{\bf\Large Are diverging CP components always nearly
proportional?}\\

\vspace{1cm} Alwin Stegeman\symbolfootnote[2]{A. Stegeman is with
the Heijmans Institute for Psychological Research, University of Groningen,
Grote Kruisstraat 2/1, 9712 TS Groningen, The Netherlands, phone: ++31 50
363 6193, fax: ++31 50 363 6304, email: a.w.stegeman@rug.nl, URL:
http://www.gmw.rug.nl/$\sim$stegeman. Research is supported by the Dutch
Organisation for Scientific Research (NWO), VIDI grant 452-08-001.} and
Lieven De Lathauwer\symbolfootnote[3]{L. De Lathauwer is with the Group
Science, Engineering and Technology of the Katholieke Universiteit Leuven
Campus Kortrijk, E. Sabbelaan 53, 8500 Kortrijk, Belgium, tel.:
+32-(0)56-32.60.62, fax: +32-(0)56-24.69.99, e-mail:
Lieven.DeLathauwer@kuleuven-kortrijk.be. He is also with the Department of
Electrical Engineering (ESAT), Research Division SCD, of the Katholieke
Universiteit Leuven, Kasteelpark Arenberg 10, B-3001 Leuven, Belgium, tel.:
+32-(0)16-32.86.51, fax: +32-(0)16-32.19.70, e-mail:
Lieven.DeLathauwer@esat.kuleuven.be, URL:
http://homes.esat.kuleuven.be/$\sim$delathau/home.html.
Research supported by: (1)
Research Council K.U.Leuven: GOA-MaNet, CoE EF/05/006
Optimization in Engineering (OPTEC), CIF1 and STRT1/08/023 (2) F.W.O.: (a)
project G.0427.10N, (b) Research Communities ICCoS, ANMMM and MLDM, (3) the
Belgian Federal Science Policy Office: IUAP P6/04 (DYSCO, ``Dynamical
systems, control and optimization'', 2007--2011), (4) EU: ERNSI.}
\\ \vspace{5mm}  \today \vspace{5mm}

\begin{abstract}
\noindent Fitting a Candecomp/Parafac (CP) decomposition (also known as
Canonical Polyadic decomposition) to a multi-way array or higher-order
tensor, is equivalent to finding a best low-rank approximation to
the multi-way array or higher-order tensor, where the rank
is defined as the outer-product rank. However, such a best low-rank
approximation may not exist due to the fact that the set
of multi-way arrays with rank at most $R$ is not closed for $R\ge 2$.
Nonexistence of a best low-rank approximation results in (groups of)
diverging rank-1 components when an attempt is made to compute the
approximation. In this note, we show that in a group of two or three
diverging components, the components converge to proportionality almost
everywhere. A partial proof of this result for larger groups of diverging
components is also given. Also, we give examples of groups of three, four,
and six non-proportional diverging components. These examples are shown to
be exceptional cases. \\~\\ \noindent{\em Keywords}: tensor decomposition,
low-rank approximation, Candecomp, Parafac, diverging components. \\~\\
\noindent {\em AMS subject classifications}: 15A18, 15A22, 15A69, 49M27,
62H25.

\end{abstract}\end{center}

\newpage
\section{Introduction}
\setcounter{equation}{0}
This note is an addendum to Stegeman \cite{S-boundary} where the following
subject is studied. Let $\circ$ denote the outer-product, and define the
rank (over the real field) of ${\cal Y}\in \R^{I\times J\times K}$ as
\be
{\rm rank}({\cal Y})=\min\{R\;|\;{\cal Y}=\sum_{r=1}^R
{\bf a}_r\circ{\bf b}_r\circ{\bf c}_r\}\,.
\ee

\noindent Let
\be
{S}_R(I,J,K)=\{{\cal Y}\in\R^{I\times J\times K}\;|\;
{\rm rank}({\cal Y})\le R\}\,,
\ee

\noindent be the set of $I\times J\times K$ arrays with at most rank $R$,
and let $\overline{S}_R(I,J,K)$ denote the closure of $S_R(I,J,K)$, i.e.
the union of the set itself and its boundary points in $\R^{I\times J\times
K}$.

Let ${\cal Z}\in\R^{I\times J\times K}$ and $||\cdot||$ denote the
Frobenius norm on $\R^{I\times J\times K}$. Consider the following low-rank
approximation problem.
\be
\label{prob-CP}
\min\{||{\cal Z}-{\cal Y}||\;|\;{\cal Y}\in{S}_R(I,J,K)\}\,.
\ee

\noindent The variables in this problem are actually the rank-1 terms in the
rank-$R$ decomposition of ${\cal Y}$:
\be
\label{CPdecomp}
{\cal Y}=\sum_{r=1}^R \omega_r\;({\bf a}_r\circ{\bf b}_r\circ{\bf c}_r)\,,
\ee

\noindent where vectors ${\bf a}_r,{\bf b}_r,{\bf c}_r$ have norm 1,
$r=1,\ldots,R$. The decomposition (\ref{CPdecomp}) is known as
Candecomp/Parafac (CP) or Canonical Polyadic decomposition. Let ${\bf
A}=[{\bf a}_1|\ldots|{\bf a}_R]$, ${\bf B}=[{\bf b}_1|\ldots|{\bf b}_R]$
and ${\bf C}=[{\bf c}_1|\ldots|{\bf c}_R]$.

Assuming rank$({\cal Z})>R$, an optimal solution of (\ref{prob-CP}) will be
a boundary point of the set ${S}_R(I,J,K)$. However, the set $S_R(I,J,K)$
is not closed for $R\ge 2$, and problem (\ref{prob-CP}) may not have an
optimal solution due to this fact; see De Silva and Lim \cite{DSL}.
This results in (groups of) diverging rank-1 terms (also known as diverging
CP components) in (\ref{CPdecomp}) when an attempt is made to compute a
best rank-$R$ approximation, see Krijnen, Dijkstra and Stegeman \cite{KDS}.
In such a case, the solution array ${\cal Y}$ converges to a boundary point
${\cal X}$ of $S_R(I,J,K)$ with rank$({\cal X})>R$. In practice this has
the following consequences: while running a CP algorithm to solve
(\ref{prob-CP}), the decrease of $||{\cal Z}-{\cal Y}||$ becomes very slow,
and some (groups of) columns of ${\bf A}$, ${\bf B}$, and ${\bf C}$ become
nearly linearly dependent, while the
corresponding weights $\omega_r$ become large in magnitude. However, the
sum of the corresponding rank-1 terms remains small and contributes to a
better CP fit. More formally, a group of diverging CP components
corresponds to an index set $D\subseteq\{1,\ldots,R\}$ such that \be
\label{eq-diverge} |\omega_r^{(n)}|\rightarrow\infty\,,\quad\quad{\rm
for\;all\;} r\in D\,, \ee \be \label{eq-small} {\rm
while}\;\;\;\;\;\|\sum_{r\in D}\omega_r^{(n)}\;({\bf a}_r^{(n)}\circ{\bf
b}_r^{(n)}\circ{\bf c}_r^{(n)})\|\quad\quad{\rm is\;bounded}\,, \ee

\noindent where the superscript $(n)$ denotes the $n$-th CP update of the
iterative CP algorithm. More than one group of diverging components may
exist. In that case (\ref{eq-diverge})-(\ref{eq-small}) hold for the
corresponding disjoint index sets.

In practice, a group of diverging components as described above is almost
always such that the corresponding columns of ${\bf A}$, ${\bf B}$, and
${\bf C}$ become nearly identical up to sign. That is, the diverging
components are nearly proportional. In this note, we focus on the question
whether this is always true or not. In Section 2, we show that in a
group of two or three diverging components, the corresponding
columns of ${\bf A}$, ${\bf B}$, and ${\bf C}$ converge to rank-1 almost
everywhere. In Section 4, we give a partial proof of this result for larger
groups of diverging components. In Sections 3, 5, and 6, we give examples of
groups of three, four, and six non-proportional diverging components,
respectively. These examples are shown to be exceptional cases.

In this note, all arrays, matrices, vectors, and scalars are real-valued.

\section{Groups of two or three diverging components}
\setcounter{equation}{0}
We need the following notation. A matrix form of the CP decomposition
(\ref{CPdecomp}) is
\be
\label{CPmatrix}
{\bf Y}_k={\bf A}\,{\bf C}_k\,{\bf B}^T\,,\quad\quad\quad
k=1,\ldots,K\,,
\ee

\noindent where ${\bf Y}_k$ is the $k$-th $I\times J$
frontal slice of ${\cal Y}$, and ${\bf C}_k$ is the diagonal matrix with
row $k$ of ${\bf C}$ as its diagonal. In (\ref{CPmatrix}), the weights
$\omega_r$ are absorbed into ${\bf A}$, ${\bf B}$, and ${\bf C}$.

We use ${\cal Y}=({\bf S},{\bf T},{\bf
U})\cdot{\cal G}$ to denote the multilinear matrix multiplication of an
array ${\cal G}\in\R^{R\times P\times Q}$ with matrices ${\bf S}$
($I\times R$), ${\bf T}$ ($J\times P$), and ${\bf U}$ ($K\times Q$).
The result of the multiplication is an $I\times J\times K$ array
${\cal Y}$ with entries
\be
y_{ijk}=\sum_{r=1}^{R}\sum_{p=1}^{P}\sum_{q=1}^{Q}
s_{ir}\,t_{jp}\,u_{kq}\,g_{rpq}\,, \ee

\noindent where $s_{ir}$, $t_{jp}$, and $u_{kq}$ are entries of ${\bf S}$,
${\bf T}$, and ${\bf U}$, respectively. We refer to multiplication $({\bf
I}_I,{\bf I}_J,{\bf U})\cdot{\cal G}$ with ${\bf U}$ nonsingular as a {\em
slicemix} of ${\cal G}$. We say that ${\cal G}$ has a {\em nonsingular
slicemix} if $I=J$ and the array $({\bf I}_I,{\bf I}_I,{\bf U})\cdot{\cal
G}$ has a nonsingular frontal $I\times I$ slice for some nonsingular ${\bf
U}$.

For later use, we state the following lemma.

\begin{lem}
\label{lem-boundarySGSD}
For $R\le\min(I,J,K)$ and ${\cal Y}\in\overline{S}_R(I,J,K)$, it holds that
${\cal  Y}=({\bf S},{\bf T},{\bf U})\cdot{\cal  G}$ for some ${\bf S}$,
${\bf T}$, ${\bf U}$ column-wise orthonormal, and some ${\cal
G}\in\overline{S}_R(R,R,R)$ with all frontal slices upper triangular.
Moreover, ${\cal  Y}\in{S}_R(I,J,K)$ if and only if ${\cal
G}\in{S}_R(R,R,R)$. \end{lem}

\noindent {\bf Proof.} See Stegeman \cite[Lemma 2.2 (b)]{S-boundary}.\ep

\noindent Our main result is the following.

\begin{thm}
\label{t-3divcomp}
Let ${\cal Y}^{(n)}=({\bf A}^{(n)},{\bf B}^{(n)},{\bf
C}^{(n)},{\bf\Omega}^{(n)})$ be an $I\times J\times K$ {\rm CP}
decomposition, i.e.
\be
\label{eq-decomp3}
{\cal Y}^{(n)}=\sum_{r=1}^R
\omega_r^{(n)}\; ({\bf a}_r^{(n)}\circ{\bf b}_r^{(n)}\circ{\bf
c}_r^{(n)})\,, \ee

\noindent with ${\bf A}^{(n)}$, ${\bf B}^{(n)}$, and ${\bf C}^{(n)}$ having
the length-$1$ vectors ${\bf a}_r^{(n)}$, ${\bf b}_r^{(n)}$, and ${\bf
c}_r^{(n)}$ as columns, respectively, and ${\bf\Omega}^{(n)}={\rm
diag}(\omega_1^{(n)},\ldots,\omega_R^{(n)})$. Let ${\cal
Y}^{(n)}\rightarrow{\cal X}$ with {\rm rank}$({\cal X})>R$, and such that
$|\omega_r^{(n)}|\rightarrow\infty$ for each $r$.
\begin{itemize}
\item[$(i)$] If $R=2$, then ${\bf A}^{(n)}$, ${\bf B}^{(n)}$, and ${\bf
C}^{(n)}$ converge to rank-$1$ matrices.

\item[$(ii)$] If $R=3$, then ${\bf A}^{(n)}$, ${\bf B}^{(n)}$, and ${\bf
C}^{(n)}$ converge to rank-$1$ matrices for almost all ${\cal X}$.
\end{itemize}
\end{thm}

\noindent {\bf Proof.} Krijnen et al. \cite{KDS} show that ${\bf A}^{(n)}$,
${\bf B}^{(n)}$, and ${\bf C}^{(n)}$ converge to rank-deficient matrices.
For $R=2$, this completes the proof. Next, let
$R=3$ and suppose $\min(I,J,K)\ge 3$. We have ${\cal
X}\in\overline{S}_3(I,J,K)$ and ${\cal X}$ is a boundary point of
$S_3(I,J,K)$. By Lemma~\ref{lem-boundarySGSD}, there exist ${\bf L}$, ${\bf
M}$, ${\bf N}$ with orthonormal columns such that ${\cal X}=({\bf L},{\bf
M},{\bf N})\cdot{\cal G}$, where the $3\times 3\times 3$ array ${\cal G}$
has all frontal slices upper triangular. We have rank$({\cal G})=$
rank$({\cal X})>3$ and the rank-3 CP sequence $({\bf L}^T,{\bf M}^T,{\bf
N}^T)\cdot{\cal Y}^{(n)}\rightarrow{\cal G}$. Slightly abusing notation, we
write ${\cal Y}^{(n)}=({\bf A}^{(n)},{\bf B}^{(n)},{\bf
C}^{(n)})\rightarrow{\cal G}$, where the weights $\omega_r^{(n)}$ have been
absorbed in ${\bf A}^{(n)},{\bf B}^{(n)},{\bf C}^{(n)}$, which are now
$3\times 3$ matrices.

We assume that ${\cal G}$ has a nonsingular slicemix. This is true for
almost all ${\cal X}$, i.e., the subset of boundary points ${\cal X}$ with
rank larger than 3, for which ${\cal G}$ does not have a nonsingular
slicemix, has lower dimensionality than the set of boundary points with
rank larger than 3. In fact, if ${\cal G}$ does not have a nonsingular
slicemix, then its upper triangular slices have a zero on their diagonals
in the same position.

We apply a slicemix to ${\cal  G}$ such that its first
slice is nonsingular. Next, we premultiply the slices of ${\cal G}$ by the
inverse of its first slice. Then ${\cal  G}$ is of the form
\be \label{eq-form3}
{\cal G}=[{\bf G}_1\,|\,{\bf G}_2\,|\,{\bf
G}_3]=\left[\begin{array}{ccc|ccc|ccc} 1 & 0 & 0 & a & d & f & \alpha &
\delta & \nu \\ 0 & 1 & 0 & 0 & b & e & 0 & \beta & \epsilon \\ 0 & 0 & 1 &
0 & 0 & c & 0 & 0 & \gamma \end{array}\right]\,. \ee

\noindent Since a matrix cannot be approximated arbitrarily well by a
matrix of lower rank, it follows that the approximating rank-3 sequence
${\cal Y}^{(n)}$ has a nonsingular slicemix for $n$ large enough. Moreover,
by Lemma~\ref{lem-boundarySGSD} we may assume without loss of generality
that ${\cal Y}^{(n)}$ has the form (\ref{eq-form3}). We denote the entries
of ${\cal Y}^{(n)}$ with subscript $n$, i.e. $a_n,\ldots,f_n$ and
$\alpha_n,\ldots,\nu_n$. Hence,
\be
\label{eq-Yn}
{\cal Y}^{(n)}=[{\bf Y}_1^{(n)}\,|\,{\bf Y}_2^{(n)}\,|\,{\bf Y}_3^{(n)}]=
\left[\begin{array}{ccc|ccc|ccc} 1 & 0 & 0 & a_n & d_n & f_n & \alpha_n &
\delta_n & \nu_n \\ 0 & 1 & 0 & 0 & b_n & e_n & 0 & \beta_n & \epsilon_n \\
0 & 0 & 1 & 0 & 0 & c_n & 0 & 0 & \gamma_n \end{array}\right]\,. \ee

\noindent Next, we consider the rank-3 decomposition $({\bf
A}^{(n)},{\bf B}^{(n)},{\bf C}^{(n)})$ of ${\cal Y}^{(n)}$, which can be
written as ${\bf Y}_k^{(n)}={\bf A}^{(n)}\,{\bf C}_k^{(n)}\,({\bf
B}^{(n)})^T$, where diagonal matrix ${\bf C}_k^{(n)}$ has row $k$ of ${\bf
C}^{(n)}$ as its diagonal, $k=1,2,3$. Since ${\bf Y}_1^{(n)}={\bf I}_3$,
matrices ${\bf A}^{(n)}$ and ${\bf B}^{(n)}$ are nonsingular. Without loss
of generality, we set ${\bf C}_1^{(n)}={\bf I}_3$. Then $({\bf
A}^{(n)})^{-1}=({\bf B}^{(n)})^T$ and ${\bf Y}_k^{(n)}={\bf A}^{(n)}\,{\bf
C}_k^{(n)}\,({\bf A}^{(n)})^{-1}$ for $k=2,3$. Hence, slices ${\bf
Y}_2^{(n)}$ and ${\bf Y}_3^{(n)}$ have the same eigenvectors. Moreover,
their three eigenvectors are linearly independent, and their eigenvalues
are on the diagonals of ${\bf C}_2^{(n)}$ and ${\bf C}_3^{(n)}$,
respectively. Since ${\bf Y}_2^{(n)}$ and ${\bf Y}_3^{(n)}$ have
eigenvalues $a_n,b_n,c_n$ and $\alpha_n,\beta_n,\gamma_n$, respectively, we
obtain
\be
\label{eq-Cn}
 {\bf
C}^{(n)}=\left[\begin{array}{ccc} 1 & 1 & 1\\ a_n & b_n & c_n \\ \alpha_n &
\beta_n & \gamma_n \end{array}\right]\,. \ee

\noindent From Krijnen et al. \cite{KDS} we know that ${\bf A}^{(n)}$,
${\bf B}^{(n)}$, and ${\bf C}^{(n)}$ converge to matrices with ranks
less than 3. The eigendecomposition ${\bf Y}_k^{(n)}={\bf A}^{(n)}\,{\bf
C}_k^{(n)}\,({\bf A}^{(n)})^{-1}$ converges to frontal slice ${\bf G}_k$
of ${\cal G}$, $k=2,3$. Hence, the eigenvectors in ${\bf A}^{(n)}$
converge to those of ${\bf G}_k$, $k=2,3$. Suppose ${\bf A}^{(n)}$ has a
rank-1 limit. Then ${\bf G}_k$ has only one eigenvector and three identical
eigenvalues, $k=2,3$. Hence, $a=b=c$ and $\alpha=\beta=\gamma$. Suppose
${\bf A}^{(n)}$ has a rank-2 limit $[{\bf a}_1\;{\bf a}_2\;{\bf a}_3]$,
with ${\bf a}_1,{\bf a}_2,{\bf a}_3$ eigenvectors associated with
eigenvalues $a,b,c$ of ${\bf G}_2$, and eigenvalues
$\alpha,\beta,\gamma$ of ${\bf G}_3$, respectively. Without loss of
generality, let ${\bf a}_1$ and ${\bf a}_2$ be linearly independent. If
${\bf a}_3$ is proportional to either ${\bf a}_1$ or ${\bf a}_2$, then
${\bf B}^{(n)}=({\bf A}^{(n)})^{-T}$ has large numbers in only two columns.
This violates the assumption of $|\omega_r^{(n)}|\rightarrow\infty$ for
$r=1,2,3$. Hence, ${\bf a}_3$ is in the linear span of $\{{\bf a}_1,{\bf
a}_2\}$ and not proportional to ${\bf a}_1$ or to ${\bf a}_2$. For an
eigenvalue $\lambda$ of ${\bf G}_k$, we define the eigenspace
\be
E_k(\lambda)=\{{\bf x}\in\R^3:\;{\bf G}_k{\bf x}=\lambda\,{\bf x}\}\,,
\quad\quad\quad k=2,3\,.
\ee

\noindent It holds that $\lambda_1\neq\lambda_2$ implies
$E_k(\lambda_1)\cap E_k(\lambda_2)=\{{\bf 0}\}$. If $a,b,c$ are distinct,
then ${\bf a}_1,{\bf a}_2,{\bf a}_3$ would be linearly independent, which
is not the case. Without loss of generality, let $a=b$. Then ${\bf a}_3\in
E_2(a)\cap E_2(c)$, which is impossible if $a\neq c$. Hence, it follows that
$a=b=c$. The proof of $\alpha=\beta=\gamma$ is analogous. This implies that
${\bf C}^{(n)}$ in (\ref{eq-Cn}) converges to a rank-1 matrix.

As ${\cal  Y}^{(n)}\rightarrow{\cal  G}$, we first assume that the
eigenvalues $a_n,b_n,c_n$ are distinct and the eigenvalues
$\alpha_n,\beta_n,\gamma_n$ are distinct. It can be verified that the
eigenvectors ${\bf A}^{(n)}$ of ${\bf Y}_2^{(n)}$ associated with
eigenvalues $a_n,b_n,c_n$ are, respectively,
\be
\label{eq-eigenvec2}
\left(\begin{array}{c} 1\\ 0\\ 0\end{array}\right)\,,\quad\quad
\left(\begin{array}{c} 1 \\ \frac{b_n-a_n}{d_n}\\ 0
\end{array}\right)\,,\quad\quad \left(\begin{array}{c} 1\\ \frac{e_n
(c_n - a_n)}{d_n e_n+f_n (c_n-b_n)}\\ \frac{(c_n-a_n)(c_n-b_n)}
{d_ne_n + f_n (c_n-b_n)}\end{array}\right)\,. \ee

\noindent As explained above, the eigenvectors of ${\bf Y}_3^{(n)}$
(in terms of $\alpha_n,\ldots,\nu_n$) must be identical to those of ${\bf
Y}_2^{(n)}$. We assume $d\neq 0$, $e\neq 0$, $f\neq 0$, $\delta\neq 0$,
$\epsilon\neq 0$, $\nu\neq 0$, which holds for almost all ${\cal X}$.
Combined with $a=b=c$ and $\alpha=\beta=\gamma$, it follows from
(\ref{eq-eigenvec2}) that ${\bf A}^{(n)}$ converges to a rank-1 matrix.
For ${\bf A}^{(n)}$ as in (\ref{eq-eigenvec2}), the columns of
${\bf B}^{(n)}=({\bf A}^{(n)})^{-T}$ are
\be
\label{eq-columnsBn}
\left(\begin{array}{c} 1\\ \frac{d_n}{(a_n-b_n)} \\
\frac{d_ne_n+f_n(a_n-b_n)}{(a_n-b_n)(a_n-c_n)}\end{array}\right)\,,
\quad\quad
\left(\begin{array}{c} 0\\ \frac{d_n}{(b_n-a_n)} \\
\frac{d_ne_n}{(a_n-b_n)(c_n-b_n)}\end{array}\right)\,,\quad\quad
\left(\begin{array}{c} 0\\ 0\\ \frac{d_ne_n+f_n (c_n-b_n)}
{(c_n-a_n)(c_n-b_n)}\end{array}\right)\,.
\ee

\noindent Hence, each column of ${\bf B}^{(n)}$ contains large numbers for
large $n$. After normalizing the third entries of each column of ${\bf
B}^{(n)}$ to 1, we obtain
\be
\label{eq-columnsBn2}
\left(\begin{array}{c} \frac{(a_n-b_n)(a_n-c_n)}{d_ne_n + f_n (a_n-b_n)}
\\ \frac{d_n (a_n - c_n)}{d_n e_n+f_n (a_n-b_n)} \\ 1\end{array}\right)\,,
\quad\quad
\left(\begin{array}{c} 0\\ \frac{(b_n-c_n)}{e_n} \\ 1\end{array}\right)\,,
\quad\quad
\left(\begin{array}{c} 0\\ 0\\ 1\end{array}\right)\,.
\ee

\noindent  It follows that also ${\bf B}^{(n)}$ converges to a rank-1
matrix.

Above, we assumed distinct eigenvalues $a_n,b_n,c_n$ and
$\alpha_n,\beta_n,\gamma_n$. Next, we show that cases with identical
eigenvalues can be left out of consideration. We only consider cases where
some of $a_n,b_n,c_n$ are identical. Cases where some of
$\alpha,\beta,\gamma$ are identical can be treated analogously.
If $a_n=b_n\neq c_n$ for $n$ large enough, then we must have $d_n=0$ to
obtain three linearly independent eigenvectors of ${\bf Y}_2^{(n)}$. This
is due to the upper triangular form of ${\bf Y}_2^{(n)}$ in (\ref{eq-Yn}).
This implies that $d=0$ in the limit, which does not hold for almost all
${\cal X}$.

The case $a_n\neq b_n=c_n$ can be dealt with analogously. Here, we must
have $e_n=0$ to obtain three linearly independent eigenvectors of
${\bf Y}_2^{(n)}$ in (\ref{eq-Yn}). This implies that
$e=0$ in the limit, which does not hold for almost all ${\cal X}$.

Next, suppose $a_n=c_n\neq b_n$ for $n$ large enough. To obtain three
linearly independent eigenvectors of ${\bf Y}_2^{(n)}$ in (\ref{eq-Yn}),
we must have $d_ne_n+f_n(c_n-b_n)=0$. Since $c_n-b_n\rightarrow
c-b=0$, this implies that $de=0$ in the limit, which does not hold for
almost all ${\cal X}$.

Finally, we consider the case $a_n=b_n=c_n$ for $n$ large enough. To
obtain three linearly independent eigenvectors of ${\bf Y}_2^{(n)}$ in
(\ref{eq-Yn}), we must have $d_n=e_n=f_n=0$. This implies that
$d=e=f=0$ in the limit, which does not hold for almost all ${\cal X}$.
This completes the proof for $\min(I,J,K)\ge 3$.

Next, let $\min(I,J,K)=2$. Without loss of generality, we assume $I\ge J\ge
K$. If $I=J=K=2$, then ${\cal X}$ is a $2\times 2\times 2$ array, which has
maximal rank 3 \cite{JJ}. A contradiction to rank$({\cal X})>3$. If
$I>J=K=2$, then by \cite[Theorem 5.2]{DSL} there exists a column-wise
orthonormal ${\bf L}$ ($I\times 3$) such that ${\cal X}=({\bf L},{\bf
I}_2,{\bf I}_2)\cdot{\cal G}$, with ${\cal G}$ a $3\times 2\times 2$ array
and rank$({\cal X})=$ rank$({\cal G})$. Since the maximal rank of $3\times
2\times 2$ arrays is 3 \cite{JJ}, we again obtain a contradiction.

Finally, let $I\ge J>K=2$. By \cite[Theorem 5.2]{DSL} there exist column-wise
orthonormal ${\bf L}$ ($I\times 3$) and ${\bf M}$ ($J\times 3$) such that
${\cal X}=({\bf L},{\bf M},{\bf I}_2)\cdot{\cal G}$, with ${\cal G}$ a
$3\times 3\times 2$ array and rank$({\cal X})=$ rank$({\cal G})>3$. The
remainder of the proof is analogous to the beginning of the proof for
$\min(I,J,K)\ge 3$. We let ${\cal Y}^{(n)}=({\bf A}^{(n)},{\bf
B}^{(n)},{\bf C}^{(n)})\rightarrow{\cal G}$, where ${\bf A}^{(n)}$ and
${\bf B}^{(n)}$ are $3\times 3$, and ${\bf C}^{(n)}$ is $2\times 3$. As
shown in \cite{Ste-arxiv}, array ${\cal G}$ can be transformed to have
upper triangular slices. Assuming ${\cal G}$ has a nonsingular slicemix, we
transform it to
\be
\label{eq-G332}
{\cal G}=[{\bf G}_1\,|\,{\bf G}_2]=\left[\begin{array}{ccc|ccc} 1 & 0 & 0 &
a & d & f \\ 0 & 1 & 0 & 0 & b & e \\ 0 & 0 & 1 & 0 & 0 & c
\end{array}\right]\,. \ee

\noindent We assume ${\cal Y}^{(n)}$ to be of the same form, with entries
$a_n,b_n,c_n,d_n,e_n,f_n$. As above, we have ${\bf B}^{(n)}=({\bf
A}^{(n)})^{-T}$,
\be
\label{eq-Cn332}
{\bf C}^{(n)}=\left[\begin{array}{ccc} 1 & 1 & 1\\ a_n & b_n & c_n
\end{array}\right]\,, \ee

\noindent and eigendecomposition ${\bf Y}_2^{(n)}={\bf A}^{(n)}\,{\bf
C}_2^{(n)}\,({\bf A}^{(n)})^{-1}$ converging to frontal slice ${\bf G}_2$.
Krijnen et al. \cite{KDS} show that ${\bf C}^{(n)}$ converges to a
rank-deficient matrix. Hence, ${\bf C}^{(n)}$ converges to a rank-1 matrix,
and we have $a=b=c$ in the limit.

We assume distinct eigenvalues $a_n,b_n,c_n$. As above, having
some identical eigenvalues for $n$ large enough yields $d=0$ or $e=0$ or
$f=0$, which does not hold for almost all ${\cal X}$. The eigenvectors
${\bf A}^{(n)}$ of ${\bf Y}_2^{(n)}$ are given by (\ref{eq-eigenvec2}). We
assume $d\neq 0$, $e\neq 0$, $f\neq 0$. Since $a=b=c$, it can be seen that
${\bf A}^{(n)}$ converges to a rank-1 matrix. The same is true for
${\bf B}^{(n)}=({\bf A}^{(n)})^{-T}$, which has columns equal to
(\ref{eq-columnsBn2}) after normalizing the third entry of each column to
1. This completes the proof. \ep

\section{Example of non-proportional diverging components for $R=3$}
\setcounter{equation}{0}
The example is of size $3\times 3\times 2$. Let
\be
{\cal X}=\left[\begin{array}{ccc|ccc} 1 & 0 & 0 & a & 0 & f \\
0 & 1 & 0 & 0 & a & e \\ 0 & 0 & 1 & 0 & 0 & a \end{array}\right]\,, \ee

\noindent with $e\neq 0$ and $f\neq 0$. Here, ${\cal X}$ plays the role of
${\cal G}$ in (\ref{eq-G332}). Since $d=0$ in ${\cal X}$ above, this is an
exception to almost all boundary arrays ${\cal X}$ with rank$({\cal X})>3$.
We have ${\cal Y}^{(n)}=({\bf A}^{(n)},{\bf B}^{(n)},{\bf
C}^{(n)})\rightarrow{\cal X}$, with
\be
\label{eq-AnBn}
{\bf A}^{(n)}=\left[\begin{array}{ccc}
1 & 0 & 1 \\
0 & 1 & \frac{e_n\,(c_n-a_n)}{f_n\,(c_n-b_n)} \\
0 & 0 & \frac{(c_n-a_n)}{f_n}\end{array}\right]\,,\quad\quad\quad
{\bf B}^{(n)}=({\bf A}^{(n)})^{-T}=\left[\begin{array}{ccc}
1 & 0 & 0 \\
0 & 1 & 0 \\
\frac{f_n}{(a_n-c_n)} & \frac{e_n}{(b_n-c_n)} & \frac{f_n}{(c_n-a_n)}
\end{array}\right]\,,
\ee

\noindent and ${\bf C}^{(n)}$ as in
(\ref{eq-Cn332}). Let $a_n=a+1/n$, $b_n=a-1/n$, $c_n=a+2/n$. Then
\be
{\bf A}^{(n)}\rightarrow\left[\begin{array}{ccc}
1 & 0 & 1 \\ 0 & 1 & \frac{e}{3f} \\ 0 & 0 & 0\end{array}\right]\,,
\quad\quad\quad
{\bf B}^{(n)}\rightarrow\left[\begin{array}{ccc}
0 & 0 & 0\\ 0 & 0 & 0\\ 1 & 1 & 1\end{array}\right]\,,
\ee

\noindent where the columns of ${\bf B}^{(n)}$ are normalized such that
their third entries are equal to 1. As we see, ${\bf A}^{(n)}$ converges to
a rank-2 matrix. Note that $|\omega_r^{(n)}|\rightarrow\infty$
for $r=1,2,3$ in this example, since all three columns of ${\bf B}^{(n)}$
in (\ref{eq-AnBn}) will have large numbers as entries.

\section{Extension to groups of four or more diverging components}
\setcounter{equation}{0}
Can statement $(ii)$ of Theorem~\ref{t-3divcomp} also be proven for $R\ge
4$ if the requirement $R\le\min(I,J,K)$ is added? A proof of this could be
analogous to the proof for $R=3$ under this requirement. That is, we have
${\cal Y}^{(n)}=({\bf A}^{(n)},{\bf B}^{(n)},{\bf C}^{(n)})\rightarrow{\cal
G}$, where ${\cal G}$ is $R\times R\times R$ and has its first slice equal
to ${\bf I}_R$ and its other slices upper triangular. Matrices ${\bf
A}^{(n)}$, ${\bf B}^{(n)}$, and ${\bf C}^{(n)}$ have size $R\times R$, with
${\bf B}^{(n)}=({\bf A}^{(n)})^{-T}$, ${\bf C}_1^{(n)}={\bf I}_R$, and
${\bf Y}_k^{(n)}={\bf A}^{(n)}\,{\bf C}_k^{(n)}\,({\bf A}^{(n)})^{-1}$
converging to frontal slice ${\bf G}_k$ of ${\cal G}$, $k=2,\ldots,R$.
Suppose we have shown that ${\bf G}_k$, $k=2,\ldots,R$, all have $R$
identical eigenvalues. Then the limit of ${\bf C}^{(n)}$ analogous to
(\ref{eq-Cn}) is a rank-1 matrix. For any fixed $R\ge 4$, we can use
symbolic computation software to obtain expressions for ${\bf A}^{(n)}$
and ${\bf B}^{(n)}$ analogous to
(\ref{eq-eigenvec2})-(\ref{eq-columnsBn2}). This would imply that also
${\bf A}^{(n)}$ and ${\bf B}^{(n)}$, when normalized to have length-1
columns, converge to rank-1 matrices. The cases where ${\bf Y}^{(n)}_k$ has
some identical eigenvalues for $n$ large enough, $2\le k\le R$, restrict
some entries of ${\cal G}$ to zero. This is an exception to almost all
boundary arrays ${\cal X}$ with rank$({\cal X})>R$.

The difficulty in obtaining the proof sketched above seems to be in showing
that ${\bf G}_k$, $k=2,\ldots,R$, all have $R$ identical eigenvalues. For
this to be proven, the possibilities for the rank and dependence structure
of the limit of ${\bf A}^{(n)}$ must be analyzed. To have a group of $R$
diverging components with $R\ge 4$, it must not only be checked that all
columns of ${\bf B}^{(n)}=({\bf A}^{(n)})^{-T}$ contain large numbers, but
also that the $R$ components do not consist of several different groups of
diverging components. For example, if $R=4$, then having two groups of two
diverging components corresponds to two times two identical eigenvalues in
the limit. In this case, the limit of ${\bf A}^{(n)}$ has two groups of two
proportional columns.

To demonstrate the above, consider the case $R=4$. It suffices to show that
${\bf G}_2$ has four identical eigenvalues. Let ${\bf G}_2$ have
eigenvalues $\lambda_1,\lambda_2,\lambda_3,\lambda_4$ and associated
eigenvectors ${\bf A}=[{\bf a}_1\;{\bf a}_2\;{\bf a}_3\;{\bf a}_4]$, with
${\bf A}^{(n)}\rightarrow{\bf A}$. Let $E(\lambda)=\{{\bf x}\in\R^4:\;{\bf
G}_2{\bf x}=\lambda\,{\bf x}\}$ denote the eigenspace corresponding to
eigenvalue $\lambda$. We have rank$({\bf A})<4$. If rank$({\bf
A})=1$, then ${\bf G}_2$ has only one eigenvector and four identical
eigenvalues: $\lambda_1=\lambda_2=\lambda_3=\lambda_4$.

Next, let rank$({\bf A})=2$. Without loss of generality, we assume ${\bf
a}_3,{\bf a}_4\in{\rm span}\{{\bf a}_1,{\bf a}_2\}$, with ${\bf a}_1$ and
${\bf a}_2$ linearly independent. Suppose $\lambda_1=\lambda_2$. Then ${\bf
a}_4\in E(\lambda_4)\cap E(\lambda_1)$, which implies
$\lambda_1=\lambda_4$. Analogously, ${\bf a}_3\in E(\lambda_3)\cap
E(\lambda_1)$ implies $\lambda_1=\lambda_3$. Hence, we obtain
$\lambda_1=\lambda_2=\lambda_3=\lambda_4$. Next, suppose
$\lambda_1\neq\lambda_2$. Because rank$({\bf A})=2$, we have at most two
distinct eigenvalues. If $\lambda_3=\lambda_1\neq\lambda_2=\lambda_4$, then
${\bf a}_1$ and ${\bf a}_3$ are proportional and ${\bf a}_2$ and ${\bf
a}_4$ are proportional. Hence, this is a case of two groups of two
diverging components, and not one group of four diverging components. If
$\lambda_1\neq\lambda_2=\lambda_3=\lambda_4$, then ${\bf a}_2,{\bf
a}_3,{\bf a}_4$ are proportional, and we have a group of three diverging
components only (i.e., large numbers in three columns of ${\bf
B}^{(n)}=({\bf A}^{(n)})^{-T}$ only). Other possibilities for
$\lambda_1\neq\lambda_2$ and rank$({\bf A})=2$ are analogous. It follows
that if rank$({\bf A})=2$, then $\lambda_1=\lambda_2=\lambda_3=\lambda_4$.

Next, let rank$({\bf A})=3$. Without loss of generality, we assume ${\bf
a}_4\in{\rm span}\{{\bf a}_1,{\bf a}_2,{\bf a}_3\}$, with ${\bf a}_1,{\bf
a}_2,{\bf a}_3$ linearly independent. Suppose
$\lambda_1=\lambda_2=\lambda_3$. Then ${\bf a}_4\in E(\lambda_4)\cap
E(\lambda_1)$, which implies $\lambda_1=\lambda_4$, and yields the desired
result. Next, suppose $\lambda_1=\lambda_2\neq\lambda_3$. If
$\lambda_4=\lambda_1$, then we have a group of three diverging components
only. If $\lambda_4=\lambda_3$, then ${\bf a}_3$ and ${\bf a}_4$ are
proportional, and we have a group of two diverging components only. If
$\lambda_4\neq\lambda_1$ and $\lambda_4\neq\lambda_3$, then rank$({\bf
A})=4$ which is not possible. Next, suppose that
$\lambda_1,\lambda_2,\lambda_3$ are distinct. Then $\lambda_4$ must be
equal to one of them. Let $\lambda_4=\lambda_1$. Then ${\bf a}_1$ and ${\bf
a}_4$ are proportional, and we have a group of two diverging components
only. Other possibilities for the equality of some eigenvalues can be
treated analogously. It follows that if rank$({\bf A})=3$, then
$\lambda_1=\lambda_2=\lambda_3=\lambda_4$.

As a final remark, we state that under the requirement $R\le\min(I,J)$, the
case $K=2$ can be proven for any $R\ge 4$ analogous to the proof of $I\ge
J>K=2$ and $R=3$.

\section{Example of non-proportional diverging components for $R=4$}
\setcounter{equation}{0}
This example is in the spirit of \cite[Section 4]{DSL}. Let $\min(I,J)\ge
6$, $K\ge 5$, and $R=4$. Let ${\bf A}=[{\bf a}_1\;{\bf a}_2]$ ($I\times 2$)
and ${\bf B}=[{\bf b}_1\;{\bf b}_2]$ ($J\times 2$) be random matrices. For
random ${\bf Q}$ ($2\times 2$), let $\tilde{\bf A}=[\tilde{\bf
a}_1\;\tilde{\bf a}_2]={\bf AQ}$ and $\tilde{\bf B}=[\tilde{\bf
b}_1\;\tilde{\bf b}_2]={\bf B}{\bf Q}^{-T}$. Then ${\bf AB}^T=\tilde{\bf
A}\tilde{\bf B}^T$, which implies
\be
\label{eq-iszero}
{\bf a}_1\circ{\bf b}_1\circ{\bf c}+{\bf a}_2\circ{\bf b}_2\circ{\bf c}-
\tilde{\bf a}_1\circ\tilde{\bf b}_1\circ{\bf c}-\tilde{\bf a}_2\circ
\tilde{\bf b}_2\circ{\bf c}={\cal O}\,,
\ee

\noindent for any vector ${\bf c}$, where ${\cal O}$ denotes an allzero
array. Let ${\bf X}=[{\bf x}_1\;\ldots\;{\bf x}_4]$ ($I\times 4$) and
${\bf Y}=[{\bf y}_1\;\ldots\;{\bf y}_4]$ ($J\times 4$) and
${\bf Z}=[{\bf z}_1\;\ldots\;{\bf z}_4]$ ($K\times 4$)
be random matrices.

We define
\be
{\bf A}^{(n)}=\left[\,
{\bf a}_1+\frac{1}{n}\,{\bf x}_1\,\right|\left.
{\bf a}_2+\frac{1}{n}\,{\bf x}_2\,\right|\left.
-\tilde{\bf a}_1-\frac{1}{n}\,{\bf x}_3\,\right|\left.
-\tilde{\bf a}_2-\frac{1}{n}\,{\bf x}_4\,\right]\,,
\ee
\be
{\bf B}^{(n)}=\left[\,
{\bf b}_1+\frac{1}{n}\,{\bf y}_1\,\right|\left.
{\bf b}_2+\frac{1}{n}\,{\bf y}_2\,\right|\left.
\tilde{\bf b}_1+\frac{1}{n}\,{\bf y}_3\,\right|\left.
\tilde{\bf b}_2+\frac{1}{n}\,{\bf y}_4\,\right]\,,
\ee
\be
{\bf C}^{(n)}=\left[\,
{\bf c}+\frac{1}{n}\,{\bf z}_1\,\right|\left.
{\bf c}+\frac{1}{n}\,{\bf z}_2\,\right|\left.
{\bf c}+\frac{1}{n}\,{\bf z}_3\,\right|\left.
{\bf c}+\frac{1}{n}\,{\bf z}_4\,\right]\,.
\ee

\noindent If we let ${\cal Y}^{(n)}=n\,({\bf A}^{(n)},{\bf B}^{(n)},{\bf
C}^{(n)})\rightarrow{\cal X}$, then by using (\ref{eq-iszero}) we obtain
\begin{eqnarray}
\label{eq-Xdecomp}
{\cal X} &=& {\bf a}_1\circ{\bf b}_1\circ{\bf z}_1+{\bf a}_1\circ{\bf
y}_1\circ{\bf c}+{\bf x}_1\circ{\bf b}_1\circ{\bf c}
+{\bf a}_2\circ{\bf b}_2\circ{\bf z}_2+{\bf a}_2\circ{\bf
y}_2\circ{\bf c}+{\bf x}_2\circ{\bf b}_2\circ{\bf c} \nonumber \\[2mm]
&& -\;\tilde{\bf a}_1\circ\tilde{\bf b}_1\circ{\bf z}_3-\tilde{\bf
a}_1\circ{\bf y}_3\circ{\bf c}-{\bf x}_3\circ\tilde{\bf b}_1\circ{\bf c}
-\tilde{\bf a}_2\circ\tilde{\bf b}_2\circ{\bf z}_4-\tilde{\bf
a}_2\circ{\bf y}_4\circ{\bf c}-{\bf x}_4\circ\tilde{\bf b}_2\circ{\bf c}\,.
\end{eqnarray}

\noindent The mode-1 rank of ${\cal X}$ equals rank$[{\bf A}\,|\,\tilde{\bf
A}\,|\,{\bf X}]=6$. The mode-2 rank of ${\cal X}$ equals rank$[{\bf
B}\,|\,\tilde{\bf B}\,|\,{\bf Y}]=6$. The mode-3 rank of ${\cal X}$ equals
rank$[{\bf c}\,|\,{\bf Z}]=5$. Since rank$({\cal X})$ is at least equal to
its mode-$i$ rank, $i=1,2,3$, it follows that rank$({\cal X})\ge 6$.
As we see, both ${\bf A}^{(n)}$ and ${\bf B}^{(n)}$ converge to rank-2
matrices, while ${\bf C}^{(n)}$ converges to a rank-1 matrix.

We create ${\cal X}$ according to (\ref{eq-Xdecomp}) and compute the
$4\times 4\times 4$ array ${\cal G}$ with upper triangular slices in the
transformation ${\cal X}=({\bf L},{\bf M},{\bf N})\cdot{\cal G}$, using
the Jacobi-type SGSD algorithm of \cite{LMV} modified as described in
\cite{SDL}. Recall that existence of this transformation follows from
Lemma~\ref{lem-boundarySGSD}. Next, we transform the slices of ${\cal
G}$ such that its first slice becomes ${\bf I}_4$. We observe the following
form for slices 2,3,4 of the transformed ${\cal G}$:
\be
\left[\begin{array}{cccc} a & 0 & e & g \\ 0 & a & c & f \\ 0 & 0 & a & 0
\\ 0 & 0 & 0 & a \end{array}\right]\,. \ee

\noindent Hence, the slices have four identical eigenvalues and zeros in
positions (1,2) and (3,4). The latter property shows that the limit point
is an exception to almost all boundary arrays ${\cal X}$ with rank$({\cal
X})>4$.

\section{Example of non-proportional diverging components for $R=6$}
\setcounter{equation}{0}
This example is similar to the one in Section 5. Let $\min(I,J,K)\ge 8$ and
$R=6$. Let
\be
{\bf A}=\left[\begin{array}{ccc} 1 & 0 & 0\\ 0 & 1 & -1\end{array}\right]\,,
\quad\quad
{\bf B}=\left[\begin{array}{ccc} 1 & 1 & 1\\ 0 & 1 & 0\end{array}\right]\,,
\quad\quad
{\bf C}=\left[\begin{array}{ccc} 1 & 0 & 0\\ 0 & 1 & 1\end{array}\right]\,,
\ee
\be
\tilde{\bf A}=\left[\begin{array}{ccc} 1 & 0 & 0\\ 1 & 1 &
-1\end{array}\right]\,, \quad\quad \tilde{\bf B}=\left[\begin{array}{ccc} 1
& 0 & 1\\ 0 & 1 & 0\end{array}\right]\,, \quad\quad \tilde{\bf
C}=\left[\begin{array}{ccc} 1 & 0 & 1\\ 0 & 1 & 0\end{array}\right]\,.
\ee

\noindent Then we have
\be
\sum_{r=1}^3 {\bf a}_r\circ{\bf b}_r\circ{\bf c}_r=\sum_{r=1}^3
\tilde{\bf a}_r\circ\tilde{\bf b}_r\circ\tilde{\bf c}_r=\left[
\begin{array}{cc|cc} 1&0&0&0\\ 0&0&0&1\end{array}\right]\,,
\ee

\noindent where the latter is a rank-2 array. Let ${\bf S}$ ($I\times
2$), ${\bf T}$ ($J\times 2$), and ${\bf U}$ ($K\times 2$) be random
matrices. Let ${\bf X}$ ($I\times 6$), ${\bf Y}$ ($J\times 6$), and ${\bf
Z}$ ($K\times 6$) be random matrices. We define
\be
{\bf A}^{(n)}=[{\bf S}\,{\bf A}\,|\,-{\bf S}\,\tilde{\bf
A}]+\frac{1}{n}\;{\bf X}\cdot{\rm diag}(1,1,1,-1,-1,-1)\,,
\ee
\be
{\bf B}^{(n)}=[{\bf T}\,{\bf B}\,|\,{\bf T}\,\tilde{\bf
B}]+\frac{1}{n}\;{\bf Y}\,,\quad\quad\quad
{\bf C}^{(n)}=[{\bf U}\,{\bf C}\,|\,{\bf U}\,\tilde{\bf
C}]+\frac{1}{n}\;{\bf Z}\,.
\ee

\noindent Since the CP decompositions $({\bf SA},{\bf TB},{\bf UC})$ and
$({\bf S}\tilde{\bf A},{\bf T}\tilde{\bf B},{\bf U}\tilde{\bf C})$ yield
the same $I\times J\times K$ array, it follows that ${\cal
Y}^{(n)}=n\,({\bf A}^{(n)},{\bf B}^{(n)},{\bf C}^{(n)})\rightarrow{\cal
X}$, with
\begin{eqnarray}
\label{eq-X6decomp}
{\cal X} &=& \sum_{r=1}^3 ( {\bf Sa}_r\circ{\bf Tb}_r\circ{\bf z}_r+
{\bf Sa}_r\circ{\bf y}_r\circ{\bf Uc}_r+{\bf x}_r\circ{\bf Tb}_r\circ{\bf
Uc}_r ) \nonumber \\
&& -\;\sum_{s=1}^3 ( {\bf S}\tilde{\bf a}_s\circ{\bf
T}\tilde{\bf b}_s\circ{\bf z}_{s+3}+ {\bf S}\tilde{\bf a}_s\circ{\bf
y}_{s+3}\circ{\bf U}\tilde{\bf c}_s+{\bf x}_{s+3}\circ{\bf
T}\tilde{\bf b}_s\circ{\bf U}\tilde{\bf c}_s )\,.
\end{eqnarray}

\noindent  The mode-1 rank of ${\cal X}$ equals rank$[{\bf
SA}\,|\,{\bf S}\tilde{\bf A}\,|\,{\bf X}]=8$. The mode-2 rank of ${\cal X}$
equals rank$[{\bf TB}\,|\,{\bf T}\tilde{\bf B}\,|\,{\bf Y}]=8$. The mode-3
rank of ${\cal X}$ equals rank$[{\bf UC}\,|\,{\bf U}\tilde{\bf C}\,|\,{\bf
Z}]=8$. Since rank$({\cal X})$ is at least equal to its mode-$i$ rank,
$i=1,2,3$, it follows that rank$({\cal X})\ge 8$. As we see, matrices ${\bf
A}^{(n)}$, ${\bf B}^{(n)}$, and ${\bf C}^{(n)}$ all converge to rank-2
matrices.

We create ${\cal X}$ according to (\ref{eq-X6decomp}) and compute the
$6\times 6\times 6$ array ${\cal G}$ with upper triangular slices in the
transformation ${\cal X}=({\bf L},{\bf M},{\bf N})\cdot{\cal G}$, using
the Jacobi-type SGSD algorithm of \cite{LMV} modified as described in
\cite{SDL}. The resulting slices of ${\cal G}$ have entries (4,4) and (5,5)
almost zero. Hence, ${\cal G}$ does not seem to have a nonsingular
slicemix. This implies that the limit point is an exception to almost all
boundary arrays ${\cal X}$ with rank$({\cal X})>6$.

\end{document}